\documentclass[a4paper,12pt]{amsart}
\usepackage{amssymb}
\usepackage{ifthen}
\usepackage{graphicx}
\usepackage{float}
\usepackage{caption}
\usepackage{subcaption}
\usepackage{cite}
\usepackage{amsfonts}
\usepackage{amscd}
\usepackage{amsxtra}
\usepackage{color}

\newtheorem{thm}{Theorem}
\newtheorem{cor}{Corollary}
\newtheorem{lem}{Lemma}

\newtheorem{conj}{Conjecture}
\newtheorem{prob}{Problem}

\theoremstyle{definition}
\newtheorem{defn}{Definition}
\newtheorem{example}{Example}
\newtheorem{rem}{Remark}

\newtheorem{case}{Case}
\newtheorem{subcase}{Subcase}

\newenvironment{pf}[1][]{%
	\vskip 1mm
	\noindent
	\ifthenelse{\equal{#1}{}}%
	{{\slshape Proof. }}%
	{{\slshape #1.} }%
}%
{\qed\bigskip}

\newcounter{alphabet}

\newenvironment{Thm}[1][]{\refstepcounter{alphabet}%
	\bigskip%
	\noindent%
	{\bf Theorem \Alph{alphabet}}%
	\ifthenelse{\equal{#1}{}}{}{ (#1)}%
	{\bf .} \itshape}{\vskip 8pt}
\newenvironment{Lem}[1][]{\refstepcounter{alphabet}%
	\bigskip%
	\noindent%
	{\bf Lemma \Alph{alphabet}}%
	{\bf .} \itshape}{\vskip 8pt}

\newcommand{\ID}{{\mathbb D}}

\def\be{\begin{equation}}
\def\ee{\end{equation}}

\newcommand{\bee}{\begin{enumerate}}
	\newcommand{\eee}{\end{enumerate}}

\newcommand{\blem}{\begin{lem}}
	\newcommand{\elem}{\end{lem}}
\newcommand{\bthm}{\begin{thm}}
	\newcommand{\ethm}{\end{thm}}
\newcommand{\bcor}{\begin{cor}}
	\newcommand{\ecor}{\end{cor}}
\newcommand{\beg}{\begin{example}}
	\newcommand{\eeg}{\end{example}}
\newcommand{\begs}{\begin{examples}}
	\newcommand{\eegs}{\end{examples}}
\newcommand{\bdefe}{\begin{defn}}
	\newcommand{\edefe}{\end{defn}}
\newcommand{\bprob}{\begin{prob}}
	\newcommand{\eprob}{\end{prob}}
\newcommand{\bques}{\begin{ques}}
	\newcommand{\eques}{\end{ques}}
\newcommand{\bei}{\begin{itemize}}
	\newcommand{\eei}{\end{itemize}}
\newcommand{\bcon}{\begin{conj}}
	\newcommand{\econ}{\end{conj}}
\newcommand{\bcons}{\begin{conjs}}
	\newcommand{\econs}{\end{conjs}}
\newcommand{\bprop}{\begin{propo}}
	\newcommand{\eprop}{\end{propo}}
\newcommand{\br}{\begin{rem}}
	\newcommand{\er}{\end{rem}}
\newcommand{\brs}{\begin{rems}}
	\newcommand{\ers}{\end{rems}}
\newcommand{\bo}{\begin{obser}}
	\newcommand{\eo}{\end{obser}}
\newcommand{\bos}{\begin{obsers}}
	\newcommand{\eos}{\end{obsers}}

\newcommand{\bca}{\begin{case}}
	\newcommand{\eca}{\end{case}}
\newcommand{\bsca}{\begin{subcase}}
	\newcommand{\esca}{\end{subcase}}

\newcommand{\bpf}{\begin{pf}}
	\newcommand{\epf}{\end{pf}}
\newcommand{\ba}{\begin{array}}
	\newcommand{\ea}{\end{array}}
\newcommand{\beq}{\begin{eqnarray}}
\newcommand{\beqq}{\begin{eqnarray*}}
\newcommand{\eeq}{\end{eqnarray}}
\newcommand{\eeqq}{\end{eqnarray*}}

\newcommand{\ds}{\displaystyle}

\def\cc{\setcounter{equation}{0}   % THIS CLEARS THE COUNTER
\setcounter{figure}{0}\setcounter{table}{0}}

\newcounter{minutes}
\divide\time by 60
\newcounter{hours}
\multiply\time by 60 \addtocounter{minutes}{-\time}
\begin{document}

\bibliographystyle{amsplain}
%
%\begin{center}
%{\tiny \texttt{FILE:~\jobname .tex,
%        printed: \number\year-\number\month-\number\day,
%        \thehours.\ifnum\theminutes<10{0}\fi\theminutes}
%}
%\end{center}

\title[Improved Bohr Inequality for harmonic mappings]{Note on Improved Bohr Inequality for harmonic mappings}

%=========================================================================
\thanks{%$^\dagger$
File:~\jobname .tex,
          printed: \number\day-\number\month-\number\year,
          \thehours.\ifnum\theminutes<10{0}\fi\theminutes}
%=========================================================================

\author[S. Ponnusamy]{Saminathan Ponnusamy
%$^\dagger $
%${}^{~\mathbf{*}}$
}
\address{
S. Ponnusamy, Department of Mathematics,
Indian Institute of Technology Madras, Chennai-600 036, India.
}
\email{samy@iitm.ac.in}

\author[R. Vijayakumar]{Ramakrishnan Vijayakumar}
\address{
R. Vijayakumar, Department of Mathematics,
Indian Institute of Technology Madras, Chennai-600 036, India.
}
\email{mathesvijay8@gmail.com}

\keywords{Bohr inequality, Bohr radius, harmonic mappings, sense-preserving  mappings,
bounded analytic functions
%\\
%$%{}^{\mathbf{*}}
%^\dagger$ {\tt This author is on leave from the Department of Mathematics,
%Indian Institute of Technology Madras, Chennai-600 036, India}
}
%\thanks{The authors thank the referees for their careful reading of the paper.
%The first author is on leave from IIT Madras, and is currently at ISI Chennai.
%}
%\date{\today  %June. 30, 09
%;  File: 2013(2).tex}

\begin{abstract}
In this paper, we give a new generalization of the Bohr inequality in refined form both for bounded analytic
functions, and for sense-preserving harmonic functions with analytic part being bounded.
\end{abstract}

\maketitle
\pagestyle{myheadings}
\markboth{S. Ponnusamy and R. Vijayakumar}{Improved Bohr Inequality for harmonic mappings}
\cc

%\subsection{Basic Notations}

\section{Introduction}
Let $\mathbb{D}=\{z\in \mathbb{C}:\,|z|<1\}$ denote the open unit disk
and $H_\infty $ denote the class of all bounded analytic functions $f$ on the unit disk $\ID$
with the supremum norm $\|f\|_\infty :=\sup_{z\in \ID}|f(z)|$. Also, let
$${\mathcal B} = \{f\in H_\infty :\, \|f\|_\infty \leq 1 \} ~\mbox{ and }~
{\mathcal B}_0=\{\omega \in {\mathcal B}:\, \omega (0)=  0 \}.
$$
Note that if $|f(z)|=1$ for some $z\in \partial\mathbb{D}$, then, by the maximum principle, it follows that $f$ should be unimodular
constant functions. So, one can conveniently exclude constant functions.
%(eg. in the discussion of subordination), since this does not affect our discussion.
 In 1914, the following theorem was proved by Harald Bohr \cite{Bohr-14} for $0\leq r\leq 1/6$ and then, it was improved to $0\leq r\leq 1/3$
independently by M. Riesz, I. Schur and F. W. Wiener.

\begin{Thm}
{\rm (Bohr, 1914)}
If $f(z)=\sum_{n=0}^{\infty} a_n z^n \in {\mathcal B}$, then
\be\label{PVW-eq0}
\sum_{n=0}^{\infty}|a_n|r^n\leq 1
\ee
holds for $0\leq r\leq 1/3$. The number $1/3$ is optimal: for $1/3<r<1$, there exists a function $f_0(z)=\sum_{n=0}^{\infty} b_n z^n \in {\mathcal B}$
such that $\sum_{n=0}^{\infty}|b_n|r^n> 1$, namely, $f_0(z) =\frac{a-z}{1-a z}$ with $a\in (0,1)$ is close to $1$ from left.
\end{Thm}

The number $1/3$ is usually referred as Bohr radius for the family ${\mathcal B}$.
Few other proofs of Bohr's theorem are known in the literature. See for instance, Paulsen, et al. \cite{PPS}, Sidon \cite{sid} and Tomic \cite{tom}.
Tomic used analogue of the argument of Sidon.
Another elementary proof of this result, which uses a result of L. Fej\'{e}r \cite{Fejer1914}, can also be found in the paper of Landau \cite{Lan1925}.
For further details, we refer to the survey articles \cite{AAP,GarMasRoss-2018} and the references therein. In the case of functions in  ${\mathcal B}_0$, the
optimal value of Bohr radius is $1/\sqrt{2}$ which was obtained by Bombieri \cite{BS1962}, where one can obtain several deep results
in this article (see also \cite{BombBour-2004}). For some refinements and investigations on Bohr, we refer to the recent articles (cf. \cite{KayPon1,KayPon3,KayPon-AASF,KayPon2}).

Let ${\mathcal F}$ consist of sequences $ \varphi =\{\varphi_{n}(r)\}_{n=0}^{\infty}$  of nonnegative continuous functions in $[0, 1)$
such that the series $\sum_{n=0}^{\infty} \varphi_{n}(r)$ converges locally uniformly with respect to $r \in [0, 1).$
%Let $\{\varphi_{n}(r)\}_{n=0}^{\infty}$ be a sequence of nonnegative continuous functions in $[0, 1)$ such that the series $\sum_{n=0}^{\infty} \varphi_{n}(r)$ converges locally uniformly with respect to $r \in [0, 1).$
Also, for $f(z)=\sum_{n=0}^{\infty} a_n z^n \in \mathcal{B}$ and $f_{0}(z):=f(z)-f(0),$ we let for convenience
$$ \Phi_N(r):=\sum_{n=N}^{\infty} \varphi_{n}(r) ~\mbox{ and }~
B_N(f, \varphi, r):= \sum_{n=N}^{\infty}|a_n| \varphi_{n}(r) \ \ \mbox{for}\ N \geq 0,
$$
so that $B_0(f, \varphi, r)=|a_0|\varphi_{0}(r)+B_1(f, \varphi, r).$ In addition, we also let
$$A(f_0,\varphi, r):=\sum_{n=1}^{\infty}|a_n|^2 \left ( \frac{\varphi_{2n}(r)}{1+|a_{0}|}+\Phi_{2n+1}(r)\right)  ~\mbox{ and }~  \|f_{0}\|^2_{r}=\sum_{n=1}^{\infty}|a_n|^2 r^{2n}.
$$
In particular, when $\varphi_n(r)=r^n$, we let
$   B_N(f, \varphi, r)=B_N(f,r)
$
and observe that  the formula for $A(f_0,\varphi, r)$ takes the following simple form
\be\label{PV-eq-ext1}
A(f_0, r):=  \left(\frac{1}{1+|a_{0}|}+\frac{r}{1-r}\right)\|f_{0}\|^2_{r} ,
\ee
since $\Phi_{2n+1}(r)=r^{2n+1}/(1-r)$.

Based on the recent investigation on this topic \cite{KaykhaPo}, the following two results were proved in \cite{PoViWirth-3}.

\begin{Thm}\label{PKV1-ThmB}
{\rm (\cite[Proof of Theorem 1]{PoViWirth-3})}
Suppose $f\in \mathcal{B},$ $f(z)=\sum_{n=0}^{\infty} a_n z^n$ and $ \{\varphi_{n}(r)\}_{n=1}^{\infty}$ belongs to ${\mathcal F}$. Then we have the following
inequality:
\beqq
B_1(f, \varphi, r)+A(f_0,\varphi, r) \leq (1-|a_0|^2)\Phi_1(r) ~\mbox{ for $r\in [0,1)$}.
\eeqq
\end{Thm}

%A sense-preserving harmonic mappings $f$ of the form $f=h+\overline{g},$ where $h$ and $g$ are analytic in $\ID$, is said to be
%$K$-quasiconformal if $|g'(z)|\leq k |h'(z)|$ in the unit disk, for $k=\frac{K-1}{K+1} \in [0,1]$.  See \cite{KayPoSha2018} for discussion on
%Bohr radius for quasiconformal harmonic mappings of the unit disk.

\begin{Lem}\label{PKV1-LemB}
%{\rm (\cite[Theorem 4 and the Proof of Theorem 5]{PoViWirth-3})}
Let $\{\psi_{n}(r)\}_{n=1}^{\infty}$ be a decreasing sequence of nonnegative   functions in $[0,r_\psi)$, and
$g, h$ be analytic in $\ID$ such that $|g'(z)|\leq k |h'(z)|$ in $\ID$ and for some $k \in [0,1],$ where $h(z)=\sum_{n=0}^\infty a_nz^n$
and $g(z)=\sum_{n=0}^\infty b_nz^n$. Then
$$\sum_{n=1}^\infty |b_n|^2 \psi_{n}(r) \leq k^2 \sum_{n=1}^\infty |a_n|^2 \psi_{n}(r) \ \ \text{for}\ \ r\in [0,r_\psi).
$$
\end{Lem}
\bpf Note that the condition $|g'(z)|\leq k |h'(z)|$  for $z\in\ID$ may be rewritten in terms of subordination as $g'(z) \prec k h'(z)$  for $z\in\ID$.
Now apply \cite[Theorem 4]{PoViWirth-3} (see also \cite[Proof of Theorem 5]{PoViWirth-3}) to obtain the desired conclusion.
\epf
%Let ${\mathcal F}_{{\rm dec}}\subset {\mathcal F}$ consist of decreasing sequences $ \{\varphi_{n}(r)\}_{n=0}^{\infty}$ of functions from ${\mathcal F}$.

In this paper, we prove two general results which in particular yield a number of recently known results as  special choices of $\varphi_{n}(r)$'s.

\section{Main Results}
\begin{thm} \label{PKV1-th1}
Suppose that $f(z)= h(z)+\overline{g(z)}=\sum_{n=0}^\infty a_nz^n +\overline{\sum_{n=1}^\infty b_nz^n}$ is harmonic mapping of the disk $\ID$
such that $|g'(z)|\leq k |h'(z)|$ in $\ID$ and for some $k \in [0,1]$. Assume that
$ \{\varphi_{n}(r)\}_{n=1}^{\infty}$ is a decreasing sequence of functions from ${\mathcal F}$,  $\Phi_1(r)=\sum_{n=1}^{\infty} \varphi_{n}(r)$ and $p\in (0,2].$ If
\beq \label{PKV1-eq2}
1 > \frac{2}{p}\left(\frac{1+r^m}{1-r^m}\right)(1+k)\Phi_1(r),
\eeq
then the following sharp inequality holds:
\be\label{PKV1-eq1}
|h(z^m)|^p+B_1(h, \varphi, r)+B_1(g, \varphi, r)+A(h_0,\varphi, r) \leq 1 \ \ \mbox{for}\ \ r \leq R^k_{m,p},
\ee
%for $r \leq R^k_{m,p},$
where $R^k_{m,p}$ is the minimal positive root of the equation
\beq \label{PKV1-eq3}
1 = \frac{2}{p}\left(\frac{1+x^m}{1-x^m}\right)(1+k)\Phi_1(x).
\eeq
In the case when $1 < \frac{2}{p}\left(\frac{1+x^m}{1-x^m}\right)(1+k)\Phi_1(x)$ in some interval $(R^k_{m,p}, R^k_{m,p}+\epsilon),$ the number $R^k_{m,p}$ cannot be improved.
\end{thm}
\begin{pf}
From the classical Schwarz-Pick lemma and Theorem~B, %\Ref{PKV1-ThmB}, 
it follows easily that
\beqq
|h(z^m)|^p+B_1(h, \varphi, r)+A(h_0,\varphi, r) \leq \left(\frac{r^m+a}{1+r^ma}\right)^p+(1-a^2)\Phi_1(r),
\eeqq
where $a=|h(0)| \in [0,1)$.
For $h \in \mathcal{B},$ as an application of Schwarz-Pick lemma, we have the
inequality $|a_{n}|\leq 1-a^2$ for $n \geq 1.$ By assumption $|g'(z)|\leq k |h'(z)|$ in $\ID,$ where $k \in[0,1]$
and so, by Lemma~C, %\Ref{PKV1-LemB}, 
it follows that
\beqq	
\sum_{n=1}^\infty |b_n|^2 \varphi_{n}(r) \leq k^2 \sum_{n=1}^\infty |a_n|^2 \varphi_{n}(r) \leq k^2 (1-a^2)^2 \Phi_1(r).
\eeqq	
Consequently, it follows from the classical Schwarz inequality that
\beqq	
\sum_{n=1}^\infty |b_n| \varphi_{n}(r) \leq  \sqrt{\sum_{n=1}^\infty |b_n|^2 \varphi_{n}(r)} \, \sqrt{\Phi_1(r)}
\leq k (1-a^2) \Phi_1(r).
\eeqq
Thus, we have
%\vspace{8pt}
%	 	
%$ \ds|h(z^m)|^p+B_1(h, \varphi, r)+B_1(g, \varphi, r)+A(h_0,\varphi, r)$
\beq \label{PKV1-eq4}
|h(z^m)|^p+B_1(h, \varphi, r)+B_1(g, \varphi, r)+A(h_0,\varphi, r)&\leq& \left(\frac{r^m+a}{1+r^ma}\right)^p+(1-a^2)(1+k)\Phi_1(r) \nonumber \\
&=&1-\Psi_{p}(a)
\eeq	
where
\beqq
\Psi_{p}(a)=1-(1-a^2)(1+k)\Phi_1(r)-\left(\frac{r^m+a}{1+r^ma}\right)^p,\ a\in [0,1].
\eeqq	
Now, as in the proof of Lemma 3.1 in \cite{LiuSamy4}, we need to determine conditions such that $\Psi_{p}(a) \geq 0$ for all $a \in [0,1].$ Note that $\Psi_{p}(1) =0.$ We claim that $\Psi_{p}$ is a decreasing function of $a,$ under the conditions of the theorem. A direct computation shows that
\beqq
\Psi'_{p}(a)=2a (1+k)\Phi_1(r)-p(1-r^{2m})\frac{(r^m+a)^{p-1}}{(1+r^m a)^{p+1}}
\eeqq	
and
\beqq
\Psi''_{p}(a)=2(1+k) \Phi_1(r)-p(1-r^{2m})\frac{(r^m+a)^{p-2}}{(1+r^ma)^{p+2}}[p-1-2ar^m-(p+1)r^{2m}].	
\eeqq
Evidently, $\Psi''_{p}(a)\geq 0$ for all $a \in [0,1],$ whenever $0<p \leq 1.$ Hence for $r \leq R^k_{m,p},$
\beqq
\Psi'_{p}(a) \leq \Psi'_{p}(1)= 2(1+k) \Phi_1(r)-p \left(\frac{1-r^m}{1+r^m}\right) \leq 0,
\eeqq
by the assumption that \eqref{PKV1-eq2} holds. Thus, for each $r \leq R^k_{m,p}$ and $0<p \leq 1,$ $\Psi_{p}$ is a decreasing function of $a \in [0,1],$
which in turn implies that $\Psi_{p}(a) \geq \Psi_{p}(1)=0$ for all $a \in [0,1]$ and the desired inequality  follows from \eqref{PKV1-eq4}.	
Next, we show that condition $\Psi'_{p}(1) \leq 0$ is also sufficient for the function $\Psi_{p}(a)$ to be decreasing on $[0,1]$
in the case when $1<p \leq 2.$
From the proof of Lemma 3.1 in \cite{LiuSamy4} it was known that
$\Phi(\sqrt[m]{r})\geq a^{p-1}$ for all $r\in [0,1)$, where
$$\Phi(r)=	(1+r^m)^2\dfrac{(r^m+a)^{p-1}}{(1+r^m a)^{p+1}},~\mbox{ for }~r \in [0,1).
$$
In view of the above discussion, for $r \leq R^k_{m,p},$ we find that
\beqq
\Psi'_{p}(a)&=& 2a(1+k) \Phi_1(r) -p\left(\frac{1-r^m}{1+r^m}\right) \Phi(r) \\
&\leq& a^{p-1}\left[2a^{2-p}(1+k)\Phi_1(r)-p\left(\frac{1-r^m}{1+r^m}\right)\right] \\
&\leq& a^{p-1}\left[2(1+k)\Phi_1(r)-p\left(\frac{1-r^m}{1+r^m}\right)\right]= a^{p-1}\Psi'_{p}(1) \leq 0,
\eeqq	
since $0 \leq a^{2-p} \leq 1$ for $1<p \leq 2.$ Again, $\Psi_{p}(a)$ is a decreasing on $[0,1],$ whenever $1<p \leq 2$ which implies that $\Psi_{p}(a) \geq \Psi_{p}(1)=0$ for all $a \in [0,1]$ and thus, the desired inequality \eqref{PKV1-eq1} holds.

Now let us prove that $R^k_{m,p}$ is an optimal number, we consider the function
$$ h(z) =\frac{a-z}{1-a z}=a-(1-a^2)\sum\limits_{n=1}^\infty a^{n-1}z^n,\quad a\in [0,1).
$$
and $g(z)= \lambda k h(z),$ where $|\lambda|=1$. Also, let $h_0(z)=h(z)-h(0)$.
Then it is a simple exercise to see that

\vspace{8pt}
	$ \ds|h(-r^m)|^p+B_1(h, \varphi, r)+B_1(g, \varphi, r)+A(h_0,\varphi, r)$
\beqq
&=& \left(\frac{r^m+a}{1+r^ma}\right)^p+(1-a^2)(1+k)\sum_{n=1}^{\infty} a^{n-1}\varphi_{n}(r) \\
&&\hspace{.2cm} +(1-a^2)^2\sum_{n=1}^{\infty} a^{2n-2} \left[ \frac{\varphi_{2n}(r)}{1+a}+\Phi_{2n+1}(r)\right]\\
&=& 1+(1-a) Q_{p}(a,r)+ O((1-a)^2),
\eeqq
where
\beqq
Q_{p}(a,r)= (1+a)(1+k)\sum_{n=1}^{\infty} a^{n-1}\varphi_{n}(r)-\frac{1}{1-a}\left(1-\left(\frac{r^m+a}{1+r^ma}\right)^p \right)  ,
\eeqq	
and it is easy to see that the last expression on the left
%$\ds |h(-r)|^p+B_1(h, \varphi, r)+B_1(g, \varphi, r)+\left(\frac{1}{1+|a_{0}|}+\Phi_1(r)\right)\|h_{0}\|^2_{r}$
is bigger than or equal to $1$ if $Q_{p}(a,r) \geq 0.$ In fact, for $r > R^k_{m,p}$ and $a$ close to $1,$ we see that
\beqq
\lim_{a \rightarrow 1^-} Q_{p}(a,r)=\left[2(1+k)\Phi_1(r)-p\left(\frac{1-r^m}{1+r^m}\right) \right]>0,
\eeqq
showing that the number $R^k_{m,p}$ in \eqref{PKV1-eq3} is best possible.
%The proof of the theorem is complete.
\end{pf}

The most fundamental special case is perhaps $\varphi _n(r)=r^n$ ($n\geq 1$) and $k=0$. In this case, Theorem \ref{PKV1-th1} gives the following.

\bcor
{\rm (\cite[Lemma 3.3]{LiuSamy4})}
Suppose that $f\in \mathcal{B}$ and $f(z)=\sum_{n=0}^{\infty}a_{n} z^{n}$ with $a=|f(0)|$ and $f_0(z)=f(z)-f(0)$. Then for $p\in (0,2]$, we have the sharp inequality:
%\begin{equation*}
%|f(z)|^p+B_1(f,r)+\left(\frac{1}{1+a}+\frac{r}{1-r}\right)\|f_0\|_r^2 \leq 1~ \mbox{ for }~ r\leq R_{p}=\frac{p}{\sqrt{4p+1}+p+1},
%\end{equation*}
%where $\|f_{0}\|^2_{r}:=\sum_{n=1}^{\infty}|a_n|^2 r^{2n}$.
\begin{equation*}
|f(z)|^p+B_1(f,r)+A(f_0, r)\leq 1~ \mbox{ for }~ r\leq R_{p}=\frac{p}{\sqrt{4p+1}+p+1},
\end{equation*}
where $A(f_0, r)$ is defined by \eqref{PV-eq-ext1}.
%where $R_{p}$ is the positive root of the equation $2(1+r)r-p(1-r)^{2}=0$.
The radius $R_{p}$ is best possible.
\ecor
%\bpf
%Set $k=0$, and  $\varphi _n(r)=r^n$ ($n\geq 1$) in Theorem \ref{PKV1-th1}.
%\epf

\begin{rem}
We mention now few other simple special cases.
\begin{itemize}
\item[(1)] For  $\varphi _n(r)=r^n$ ($n\geq 1$),   Theorem \ref{PKV1-th1}
gives a refinement of \cite[Theorem 5]{AlkhPo-2020}.
%\item[(1)] If we set $m=1$ in Corollary \ref{cor 1}, then we easily have the sharp inequality
%\beqq
%|h(z)|^p+B_1(h, \varphi, r)+B_1(g, \varphi, r) \leq 1\ \mbox{for}\ r \leq R^k_{1,p},
%\eeqq
%where $0<p\leq 2$ and $R^k_{1,p}$ is the minimal positive root of the equation
%\beqq
%1 = \frac{2}{p}\left(\frac{1+r}{1-r}\right)(1+k)\Phi_1(x).
%\eeqq
%\item[(2)] For $\varphi_{n}(r)=r^n \, (n \geq 1)$ and
In this special case, if we allow $m \rightarrow \infty$, then the resulting inequality \eqref{PKV1-eq1} takes the form
%$R^k_{m,p} \rightarrow R(p):=\dfrac{p}{2(1+k)+p}.$ Thus, under the hypotheses of Corollary \ref{cor 1}, we easily have
\beq \label{PKV1-eq5}
|h(0)|^p+\sum_{n=1}^\infty |a_n|r^n+\sum_{n=1}^\infty |b_n|r^n +A(h_0, r) \leq 1
\eeq
for $r \le R^k_{p}=\dfrac{p}{2(1+k)+p}$. The constant $R^k_{p}$ cannot be improved. In particular, $p=1,2$ in \eqref{PKV1-eq5} with the substitution
$k=\frac{K-1}{K+1}$, the inequality yields the following refined forms of two well-known results, namely,
\cite[Theorem 1.1]{KayPoSha2018} and \cite[Theorem 1.2]{KayPoSha2018}, respectively:
\begin{itemize}
\item[(a)] $\ds	|h(0)|+\sum_{n=1}^\infty |a_n|r^n+\sum_{n=1}^\infty |b_n|r^n +A(h_0, r) \leq 1 \  \ \text{for}\  \ r \le \frac{K+1}{5K+1}.$
The constant $\frac{K+1}{5K+1}$ is sharp.
		
\item[(b)] $\ds |h(0)|^2+\sum_{n=1}^\infty |a_n|r^n+\sum_{n=1}^\infty |b_n|r^n +A(h_0, r) \leq 1 \  \ \text{for}\  \ r \le \frac{K+1}{3K+1}.$
The constant $\frac{K+1}{3K+1}$ is sharp.
\end{itemize}
\item[(2)]  The case $k=0$ of \eqref{PKV1-eq5} gives the refined form of the classical Bohr's inequality and the case $k=1,$ contains the refined
Bohr inequality for sense-preserving harmonic mappings  $f(z)= h(z)+\overline{g(z)}$
of the disk $\mathbb{D}$ (see \cite[Corollary 1.4]{KayPoSha2018}).	
\end{itemize}
\end{rem}	

In the case of analytic functions, one can easily get the following result and,
since the proof of it follows on the similar lines of the proof of Theorem \ref{PKV1-th1},
we omit the details.	

\begin{thm}\label{PKV-th2}%{thm 2}
Let $f\in \mathcal{B},$ $f(z)=\sum_{n=0}^{\infty} a_n z^n$, $0<p\leq2$ and $ \{\varphi_{n}(r)\}_{n=0}^{\infty}$ belongs to ${\mathcal F}$. If
\beqq
1 > \frac{2}{p}\left(\frac{1+r^m}{1-r^m}\right)\Phi_N(r),
\eeqq
for some $N \geq 1$ then the following sharp inequality holds:
\beqq
|f(z^m)|^p+B_N(f, \varphi, r) \leq 1\ \ \mbox{for all}\ \ r \leq \rho^N_{m,p},
\eeqq
where  $\rho^N_{m,p}$ is the minimal positive root of the equation
\beqq
1 = \frac{2}{p}\left(\frac{1+x^m}{1-x^m}\right)\Phi_N(x).
\eeqq
In the case when $ 1 < \frac{2}{p}\left(\frac{1+x^m}{1-x^m}\right)\Phi_N(x)$ in some interval $(\rho^N_{m,p}, \rho^N_{m,p}+\epsilon),$ the number $\rho^N_{m,p}$ cannot be improved.	
\end{thm}
		
Note that for $N=1$, the number $\rho^N_{m,p}$ coincides with the number $R^0_{m,p}$ of Theorem \ref{PKV1-th1}.
\begin{rem}
We mention now several useful remarks concerning some special cases of Theorem \ref{PKV-th2} for $f \in \mathcal{B}$ with $f(z)= \sum_{n=0}^{\infty} a_{n}z^{n}$.
\begin{itemize}
\item[(1)] For $\varphi_{n}(r)=r^n \, (n \geq 1)$, we easily have the following sharp inequality
\beq \label{eq 5}
|f(z^m)|^p+ \sum_{n=N}^{\infty} |a_n| r^n \leq 1\ \ \mbox{for all}\ \ r \leq \rho^N_{m,p},
\eeq
where $0<p\leq2$ and $\rho^N_{m,p}$ is the minimal positive root of the equation
\beqq
2(1+r^m)r^N - p(1-r)(1-r^m) =0.
\eeqq
In particular, the case $m=1$ and $p=1,2$ in \eqref{eq 5} yields the well-known result (see \cite[Theorem 1]{KaySamy}).

\item[(2)] For $p=1$ and $\varphi_{n}(r)=r^n \, (n \geq 1)$, we obtain the result \cite[Theorem 2]{KaySamy}.

\item[(3)] For $N=1$, and $\varphi_{n}(r)=r^n \, (n \geq 1)$,  if we allow $m \rightarrow \infty$ in Theorem \ref{PKV-th2}, we obtain that
$\rho ^1_{m,p} \rightarrow \dfrac{p}{2+p}.$ Thus,  we easily have the following inequality which contains  the classical Bohr
inequality (i.e., the case $p=1$):
\beqq
|a_0|^p+ \sum_{n=1}^\infty |a_n|r^n \leq1\ \mbox{for} \ r \leq R(p)=\dfrac{p}{2+p}
\eeqq
and the constant $\dfrac{p}{2+p}$ cannot be improved.

\item[(4)] For the case $m=1=p$,  $\varphi_{2n}(r)=r^{2n} \, (n \geq 1)$ and $\varphi_{2n-1}=0 \, (n \geq 1)$, we easily have
\beqq
|f(z)|+ \sum_{n=1}^\infty |a_{2n}|r^{2n} \leq 1 ~\mbox{ for  }~ r \leq R=\sqrt{2}-1
\eeqq
and the radius $R=\sqrt{2}-1$ is the best possible. This was obtained in \cite[Lemma 2.8]{LiuShangXu}.

\item[(5)] If we allow $m \rightarrow \infty$, with $\varphi_{kn}(r)=r^{kn} \, (n \geq 1)$ for each fixed $k\geq 1$ and $\varphi_{m}=0$ for $m\neq kn$, we obtain
the following
%: if $f(z)=\sum_{n=0}^\infty a_{kn} z^{kn}$ is analytic in $\mathbb{D}$, and $\left|f(z)\right|\le 1$ in $\mathbb{D}$, then we have
the sharp inequality
\beqq
|a_0|^p+ \sum_{n=1}^\infty |a_{kn}|r^{kn}\leq1\ \mbox{for} \ r \leq \sqrt[k]{\dfrac{p}{2+p}}.
\eeqq
For $p=1$, this gives \cite[Lemma 2.1]{AliBarSoly}. For $p=2$ and $k=1$, this is a well-known result (cf. \cite[Corollary 2.7]{PPS}).
\end{itemize}

\end{rem}

\end{document}